\documentclass[12pt]{article}
\usepackage{latexsym}
\usepackage{amsfonts}
\usepackage{epsfig}
\usepackage{graphicx}

\newcommand{\ba}{\begin{array}}\newcommand{\ea}{\end{array}}
\newcommand{\sym}{S}
\newcommand{\ns}{\rm}

\renewcommand{\Bbb}{\mathbb}
\newcommand{\nse}{\kern-3pt\ns$=$}\newcommand{\qd}{\hfill$\Box$\medbreak}

\newcommand{\Z}{\mbox{$\Bbb Z$}}
\newcommand{\C}{\mbox{$\Bbb C$}}
\newcommand{\noi}{\noindent}
\newcommand{\ext}{\raise1pt\hbox{$\ts\bigwedge$}}
\renewcommand{\sym}{S}
\newcommand{\ts}{\textstyle}
\newcommand{\rf}[1]{(\ref{#1})}\renewcommand{\c}[1]{\cite{#1}}
\newcommand{\chii}{\raise2pt\hbox{$\chi$}}

\newcommand{\Fg}{\mbox{${\cal F}\kern-2pt_g$}}

\newcommand{\Mg}{\mbox{${\cal M}\kern-2pt_g$}}
\newcommand{\Ng}{\mbox{${\cal N}\kern-2pt_g$}}
\newcommand{\V}{V\kern-1pt}

\newcommand{\Gg}{\mbox{${\cal G}\kern-2pt_g$}}
\newcommand{\cir}{\raise1.6pt\hbox{\footnotesize$\circ$}}

\newcommand{\ind}{\mbox{\ns ind}}

\newcommand{\mod}{\mbox{\ns mod\,\,\,}}

\newcommand{\Res}[2]{\hbox{\ns Res}\kern-16pt\lower5pt\hbox{\footnotesize$_{#1}$}\kern2pt\left[#2\right]}
\newcommand{\qk}{quaternion-K\"ahler\kern2pt}\renewcommand{\,}{\kern1pt}
\newcommand{\ch}{\mbox{\ns ch}}

\newcommand{\coker}{\mbox{\ns coker}}
\newcommand{\td}{\hbox{Td}}

\newcommand{\dirac}{/\kern-5pt\partial}

\newcommand{\lra}{\longrightarrow}
\renewcommand{\ts}{\textstyle}
0
\newtheorem{theo}{Theorem}[section]\newcommand{\A}{\mbox{$\widehat A$}}
\newtheorem{defi}{Definition}[section]
\newtheorem{remark}{Remark}[section]
\newcommand{\Ol}{\mbox{${\cal O}$}}
\newtheorem{prop}{Proposition}[section]\def\frac#1#2{{#1\over#2}}
\def\be#1\ee{\begin{equation}#1\end{equation}}

\begin{document}
\title{Complex contact manifolds and $S^1$ actions}

\author{
Hayde\'e Herrera\footnote{Department of Mathematical Sciences,
Rutgers University, Camden, NJ 08102, USA.\hfill\break E-mail:
haydeeh@camden.rutgers.edu}  \,\,\,\,\, and
Rafael
Herrera\footnote{Centro de Investigaci\'on en Matem\'aticas, A. P.
402, Guanajuato, Gto., C.P. 36000, M\'exico. E-mail:
rherrera@cimat.mx} \footnote{Partially supported by
a UCMEXUS grant, and
CONACYT grants:
J48320-F, J110.387/2006.
} }

\date{}

\maketitle


{
\abstract{

We prove rigidity and vanishing theorems for several holomorphic Euler
characteristics on complex contact manifolds
admitting holomorphic circle
actions preserving the contact structure.
Such vanishings are reminiscent of those of LeBrun and Salamon on
Fano contact manifolds but under a symmetry assumption instead of
a curvature condition.

}
}

\section{Introduction}

The geometry of complex contact manifolds was first studied by
Kobayashi \c{Kobayashi} and Boothby \c{Boothby}, and more recently
by LeBrun \c{LeBrun-contact} and Moroianu \c{Moroianu-contact} in relation
to quaternion-K\"ahler geometry.
Here, we study these manifolds from the point of
view of transformation groups.

Inspired by Atiyah and Hirzebruch \c{AH}, Hattori proved the
vanishing of indices of Dirac operators with coefficients in
certain powers of the Spin$^c$ complex line bundle on compact
Spin$^c$ manifolds admitting smooth circle actions \c{Hattori}.
Such vanishings apply to complex contact manifolds since their
first Chern class is a multiple of an integral cohomology class.

In this note, we prove the vanishing of several holomorphic Euler
characteristics on complex contact manifolds admitting holomorphic
circle actions preserving the contact structure. The vector
bundles considered in the holomorphic Euler characteristics are
tensor products of a suitable exterior power of the contact
(distribution) sub-bundle
and a power of the canonical line bundle.

These vanishings are reminiscent of those of LeBrun and Salamon for Fano
contact manifolds \c{LS}. Their vanishings depend on a positive-curvature
condition (Fano condition) which, in particular, makes such manifolds projective.
Here, we assume the existence of a compatible circle action on a complex
contact manifold; in principle, the manifolds may neither fulfill a curvature
condition nor be projective.

The note is organized as follows.
In Section~\ref{complex-contact-manifolds}, we recall the
definition and some properties of complex contact manifolds,
the rigidity of elliptic operators,
state our main theorem (see Theorem~\ref{main-theorem}), and
describe some properties of
the exponents of the action.
In Section~\ref{index-calculations}, we carry out index
calculations and prove
Theorem~\ref{main-theorem}.
In Section~\ref{special-case} we prove further vanishings under
a non-negativity assumption on the exponents of the action.

\vspace{.1in}

\noi {\em Acknowledgements}.  The second named
author wishes to thank
the Max Planck Institute for Mathematics (Bonn) and the
Centre de Recerca Matematica (Barcelona) for their hospitality and
support.

\section{Preliminaries}\label{complex-contact-manifolds}

\subsection{Complex contact manifolds}
Let $X$ be a complex manifold and $TX$ its holomorphic tangent
bundle.

\begin{defi}
The { complex manifold} $X$ is called {\em contact} if there is a
 complex-codimension one
holomorphic sub-bundle $D$ of $TX$ which is maximally
non-integrable, i.e. the tensor
\begin{eqnarray}
D\times D &\lra& TX/D\nonumber\\ (v,w) &\mapsto& [v,w] \,\,\,\mod
D \nonumber
\end{eqnarray}
is  non-degenerate for every point of $X$.
\end{defi}

\vspace{.1in}

\noi {\bf Examples}.
\begin{itemize}
\item Let $V$ be a compact
complex manifold. Then the projectivization of the cotangent
bundle $X:= \mathbb{P}(T^*V)$ is a complex contact manifold (see
\cite{Kobayashi} for further details). Here a 1-form $\theta$ can be
defined as follows: $\theta (u) := v(d\pi(u))$, where $u\in
T_v(T^*V)$ and $\pi: T^*V\rightarrow V$ is the projection onto
$V$. Thus, $D:= ker(\theta)$.

\item Let $M$ be a quaternion-K\"ahler manifold. Its twistor space
$Z$ is a contact complex manifold (see \cite{LS}). In fact, $Z$ is
a fiber bundle over $M$ with fiber $\mathbb{C}P^1$, and $D$ is a
complex codimension 1 distribution that is transverse to the
fibers of $Z\rightarrow M$.
\end{itemize}

Let $L:= TX/D$ be the quotient line bundle and $\theta : TX \lra
L$ the tautological projection, so that we have the short exact
sequence
\begin{equation}
0\lra D \lra TX \lra L \lra 0.\label{D-TX-L}
\end{equation}

The projection $\theta$ can be thought of as a 1-form with values
in the line bundle $L$, $\theta\in\Gamma(X, \Omega^1(L))$, with
$\ker(\theta)=D$. The sub-bundle $D$ must have even rank $2n$ and,
therefore, the manifold $X$ has odd complex dimension $2n+1\geq
3$. Moreover, the non-degeneracy condition implies
\[\theta\wedge (d\theta)^n \in \Gamma (X,\Omega^{2n+1}(L^{n+1}))\]
is nowhere zero. This provides an isomorphism \c{Kobayashi,LS} of
the anti-canonical line bundle of $X$ and $L^{n+1}$,
\[\kappa_X^{-1}=\ext^{2n+1}TX \cong L^{n+1}.\]
Since $L=TX/D$, there is a $C^{\infty}$ isomorphism
\[
TX \cong D \oplus L,
\]
so that
\[c(X)=c(D)\cdot c(L).\]
There is also the following isomorphism (cf. \cite[p. 116]{LS})
\begin{eqnarray}
D\cong D^*\otimes L.\label{magic}
\end{eqnarray}

By means of the splitting principle we can write the Chern classes
in terms of formal roots
\[c(D)=(1+y_1)(1+y_2)\cdots(1+y_{2n}),\]
and
\[c(L)=(1+y_{2n+1}),\]
so that
\[c_1(X)=(n+1)y_{2n+1}.\]

\subsection{Rigidity of elliptic operators}
Let $M$ be a compact manifold and $E$ and $F$ vector bundles over
$M$.
\begin{defi}
Let $\mathcal{D}: \Gamma(E) \lra \Gamma(F)$ be an elliptic
operator acting on sections of  $E$ and $F$. The index of
$\mathcal{D}$ is the virtual vector space $\ind
(\mathcal{D})=\ker(\mathcal{D})-\coker(\mathcal{D})$. If $M$
admits a circle action preserving $\mathcal{D}$, i.e. such that
$S^1$ acts on $E$ and $F$, and commutes with $\mathcal{D}$,
$\ind(\mathcal{D})$ admits a Fourier decomposition into complex
$1$-dimensional irreducible representations of $S^1$
$\ind(\mathcal{D}) =\sum a_m\,\, L^m$, where $a_m\in \Z$ and $L^m$
is the representation of $S^1$ on $\C$ given by $z \mapsto z^m$.
The elliptic operator $\mathcal{D}$ is called {\em rigid} if
$a_m=0$ for all $m\neq 0$, i.e. $\ind(\mathcal{D})$ consists only
of the trivial representation with multiplicity $a_0$.
\end{defi}

\begin{remark}{\rm Equivalently, we can take the trace for $z\in
S^1$,
\[ \ind(\mathcal{D})_z= \mbox{trace}(z, \sum a_m\,\, L^m)=\sum a_m\,\, z^m
,\]to get a finite Laurent series in $z$. Now $\mathcal{D}$ is
rigid if and only if $\ind(\mathcal{D})_z$ does not depend on
$z\in S^1$. }\end{remark}

\noi {\bf Example}. The deRham complex
\[ \mathcal{D}=d+d^* \colon \Omega^{even} \lra \Omega^{odd}\]
from even-dimensional forms to odd-dimensional ones, where $d^*$
denotes the adjoint of the exterior derivative $d$,  is {\em rigid}
for any circle action on $M$ by isometries since by Hodge theory
the kernel and the cokernel of this operator consist of harmonic
forms, which by homotopy invariance stay fixed under the circle
action.

\subsection{Rigidity and vanishing theorem}

Now we can state the main theorem,

\begin{theo}\label{main-theorem}
Let $X$ be a complex contact manifold, $D$ the contact
distribution and $L=TX/D$. Assume $X$ admits a circle
action by holomorphic automorphisms preserving the contact
structure.
Then, the equivariat holomorphic Euler characteristic
$\chii(X,\Ol(\ext^pD^*\otimes L^{-k}))_z$ is rigid, i.e.
\[\chii(X,\Ol(\ext^pD^*\otimes L^{-k}))_z=\chii(X,\Ol(\ext^pD^*\otimes L^{-k}))\]
for all $z\in S^1$, if
\begin{eqnarray}
&0\leq k\leq n+1-p,&\quad \hfill\mbox{for $0\leq p\leq
n$,}\nonumber\\ &n-p\leq k\leq 1,&\quad \hfill\mbox{for $n+1\leq
p\leq 2n$.}\nonumber
\end{eqnarray}
Furthermore
\[\chii(X,\Ol(\ext^pD^*\otimes L^{-k}))=0\]
if
either side of the corresponding inequality is strict.
\end{theo}

\vspace{.1in}

We postpone the proof of the theorem until
Section~\ref{index-calculations} while we recall other
preliminaries.

\begin{remark}{\rm
The bounded and dotted region in the $(k,p)$-plane in Figure~\ref{fig:VanishingRegion},
shows the pairs of powers that give rigidity and
vanishing theorems for the holomorphic Euler characteristics
$\chii(X,\Ol(\ext^p D^*\otimes L^{-k}))$
when $n=5$.

\begin{figure}[ht]
\centering
\includegraphics[height=80mm]{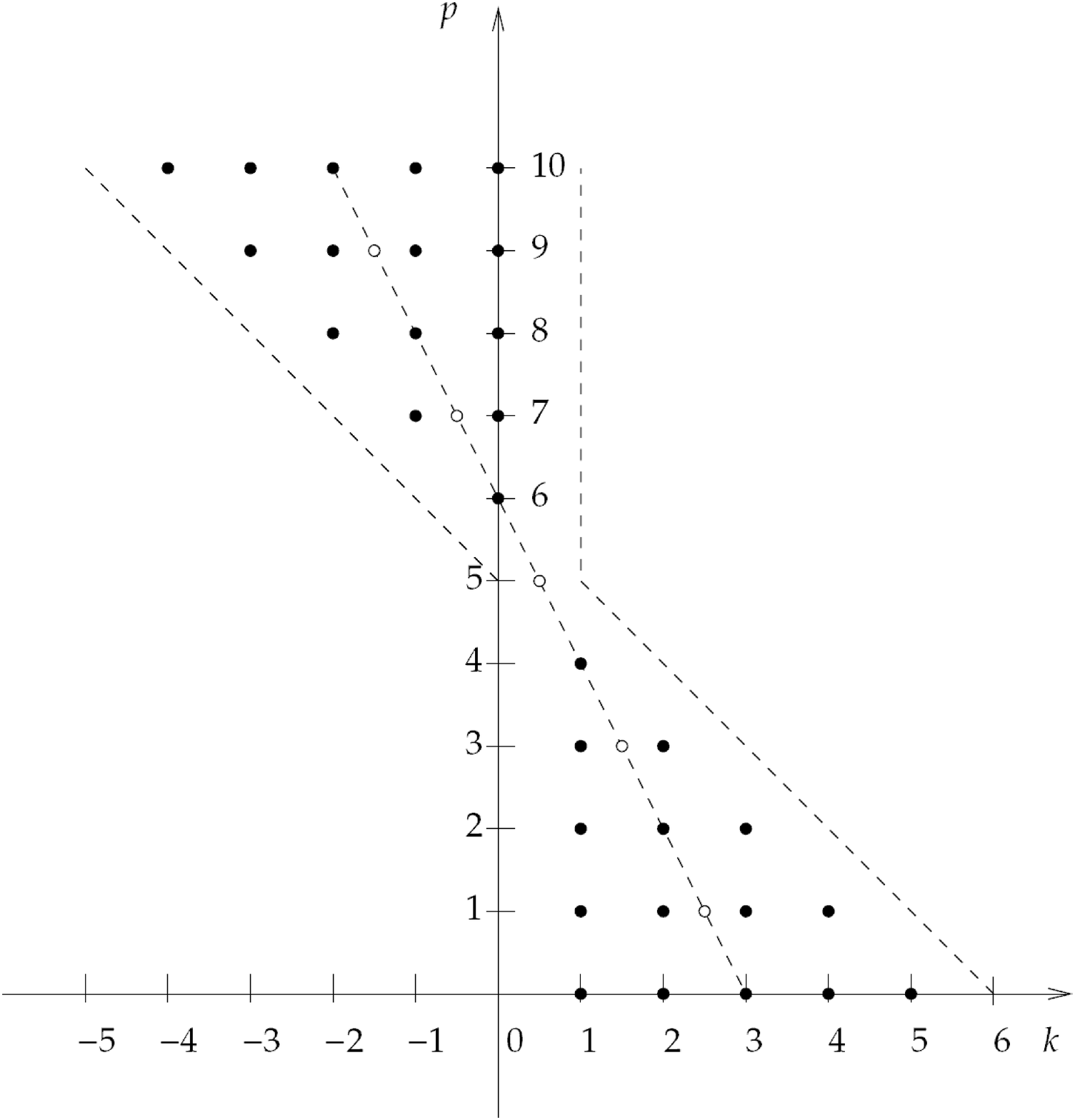}
\caption{Rigidity and vanishing region for $n=5$.} \label{fig:VanishingRegion}
\end{figure}



}
\end{remark}


\subsection{Exponents of the circle action}\label{fixed-points}

From now on, we shall assume that the complex contact manifold $X$
admits a circle action by holomorphic automorphisms preserving the
contact distribution.

Let $N$ denote a connected
component of the $S^1$-fixed point set $M^{S^1}$, which is a submanifold
and has even real codimension. Let $x\in N$. Since the
$S^1$-action preserves the contact structure, we have
$S^1$-representations on the complex vector spaces $T_xX$, $D_x$
and $L_x$ given by the fibers of the bundles $TX$, $D$ and $L$ at
$x$, as well as an $S^1$-equivariant exact sequence
\[ 0\lra D_x \lra T_xX \lra L_x \lra 0.\]
First, let us consider $T_xX$. It decomposes as a finite direct sum of
$S^1$-representations
\begin{equation}
T_xX = \bigoplus_{m\in\mathbb{Z}}
V(m).\label{decomposition-TX}
\end{equation}
where for each $m\in\mathbb{Z}$, $v\in V(m)$ and $z\in S^1$,
$z$ acts on $v$ by multiplication with $z^m$.
Similarly,
$D_x\subset T_xX$ will consist of some of these summands.
Finally, since
$\ext^{2n+1}TX \cong L^{n+1}$,
\[L_x = \bigotimes_{m} \ext^{\dim V(m)} V(m).\]
The exponents $m$ depend on the connected component $N$.
In order to carry out our computations, we will consider each
$V(m)$ to be the sum of appropriately chosen one dimensional
representations of $S^1$ with the same exponent $m$, and will
make no reference to $V(m)$ anymore.

Thus, the holomorphic tangent bundle of $X$ restricted to $N$ splits as
a sum of $S^1$-equivariant line bundles. We can think of
such a splitting as follows:
\[TX|_N ={\cal L}^{m_1}\oplus\ldots\oplus {\cal L}^{m_{2n+1}},\]
where $m_i\in \mathbb{Z}$, ${\cal L}^{m_i}$ denotes the line
bundle whose fiber is acted on by $z\in S^1$ by multiplication
with $z^{m_i}$.
Furthermore, the
lines with exponent equal to 0 add up to the tangent bundle of $N$.

Let $x\in N$. Since $\ext^{2n+1} TX = L^{n+1}$, $L_x$ has exponent
\begin{equation}
h={1\over n+1}\left( \sum_{i=1}^{2n+1} m_i\,\,\,\right).\label{h-con-TX}
\end{equation}
Since the exponents of $D_x$ are $m_1,\ldots, m_{2n}$, the exponents
of $D_x^*$ must be $-m_1,\ldots, -m_{2n}$. By \rf{magic}, the
exponents of $D_x^* = D_x\otimes L_x^{-1} $ are
\begin{eqnarray}
m_1-h &=& -m_{\sigma(1)}\nonumber\\
 &\vdots&\nonumber\\
m_{2n}-h &=& -m_{\sigma(2n)},\nonumber
\end{eqnarray}
where $\sigma$ is a permutation of $\{1, 2, \ldots, 2n\}$,
depending on the connected component $N$.
The relevant fact here is
\begin{equation}
m_i + m_{\sigma(i)} = h.\label{main-identity}
\end{equation}
Furthermore,
\[\ext^{2n}D^* = \ext^{2n}D\otimes L^{2n}\]
which implies
\[2nh = 2(m_1+m_2+\ldots + m_{2n}),\]
i.e.
\begin{equation}
nh = m_1+m_2+\ldots + m_{2n}.\label{h-con-D}
\end{equation}
Combining \rf{h-con-TX} with \rf{h-con-D}
\begin{equation}
h = m_{2n+1}.\label{h=m2n+1}
\end{equation}

\section{Proof of Theorem~\ref{main-theorem}}\label{index-calculations}

We will use the notation set up in the previous section.
Let us now consider the Hilbert polynomials in two variables given by
the following holomorphic Euler characteristics
\[\chii(X,\Ol(\ext^p D^*\otimes L^{-k}  )),\]
where $p,k\in\mathbb{Z}$, and $0\leq p\leq 2n$.
By the Atiyah-Singer index theorem, they can be computed by the following
{\small
\[\big< \ch(\ext^p D^*\otimes L^{-k})\,\,\td(X),[X]
\big>=\hspace{3in}
\]
\begin{eqnarray}
&=&\big< \ch(\ext^p
D^*)\,e^{-ky_{2n+1}}\,e^{(n+1)y_{2n+1}/2}\,\A(X),[X]\big>\nonumber\\
&=&\big< \ch(\ext^p
D^*)\,e^{(-2k+n+1)y_{2n+1}/2}\,\prod_{i=1}^{2n+1}
{y_i/2\over\sinh(y_i/2)},[X]  \big>\nonumber\\ &=&\left<
\left(\sum_{1\leq i_1<\ldots<i_p\leq 2n}
e^{-y_{i_1}-\ldots-y_{i_p}}\right)
\,e^{(-2k+n+1)y_{2n+1}/2}\,\prod_{i=1}^{2n+1}
{y_i/2\over\sinh(y_i/2)},[X]  \right>,\label{topological-index}
\end{eqnarray}
}

\noindent where $\ch$, $\td$ and $\A$ denote the Chern character, the Todd genus and the
$\A$-genus respectively.

Since the manifold $X$ admits a holomorphic $S^1$ action
preserving the contact structure, we can apply the
Atiyah-Singer fixed point theorem  \c{AS3} to obtain the equivariant
version of the index for $z\in S^1$ (cf. \cite[p. 67]{HiBJ})
{\small
\[\chii(X,\Ol(\ext^p D^*\otimes L^{-k}  ))_z=\hspace{3in}\]
\begin{eqnarray}
{}&=&\sum_{N\in X^{S^1}}\left< \left(\sum_{1\leq i_1<\ldots<i_p\leq
2n}
z^{m_{i_1}+\ldots+m_{i_p}}e^{-y_{i_1}-\ldots-y_{i_p}}\right)\,\,\times\right.\nonumber\\
&&\left.\times \,\,\, z^{(2k-n-1)h/2}e^{(-2k+n+1)y_{2n+1}/2}\,
\left(\prod_{m_i=0} {y_i\over 2}\right)
\left(\prod_{i=1}^{2n+1}
{1\over z^{-m_i/2}e^{y_i/2}-z^{m_i/2}e^{-y_i/2}}\right),[N]  \right> \nonumber\\
&=&\sum_{N\in X^{S^1}}\left< \left(\sum_{1\leq i_1<\ldots<i_p\leq
2n}
e^{-y_{i_1}-\ldots-y_{i_p}}\cdot e^{(-2k+n+1)y_{2n+1}/2}\,
\left(\prod_{m_i=0} {y_i\over 2}\right)  \times    \right.\right.\nonumber\\
&&\left.\times\left. \,\,\,
z^{m_{i_1}+\ldots+m_{i_p}+kh}
\left(\prod_{i=1}^{2n+1}
{1\over e^{y_i/2}-z^{m_i}e^{-y_i/2}}\right)\right),[N]  \right>, \label{localization-formula}
\end{eqnarray}
}

\noindent where we have substituted
\begin{eqnarray}
M &\mbox{by}& N,\nonumber\\
e^{\pm y_i} &\mbox{by}& z^{\mp m_i}e^{\pm y_i}.
\nonumber
\end{eqnarray}
in the formula \rf{topological-index} for
$\chii(X,\Ol(\ext^p D^*\otimes L^{-k}  ))$
in order to obtain the formula \rf{localization-formula}.

\begin{remark}\label{the-remark}
{\rm We wish to control the behaviour of \rf{localization-formula}
at $0$ and $\infty$ thought of as a rational function in $z$. Let
us consider the function of $z\in \mathbb{C}$
\[F(z)={z^{l}\over e^{x} - z^{m}e^{-x}}\]
where $x$ is an unknown and $l,m\in\mathbb{Z}$. First, let us assume
$m> 0$.
Thus, if $l>0$ then
\begin{equation}
\lim_{z\rightarrow 0 } F(z) = \lim_{z\rightarrow 0} {z^{l}\over
e^{x} - z^{m}e^{-x}} =0,   \label{zero1}
\end{equation}
 and if $l-m<0$ then
\begin{equation}
\lim_{z\rightarrow \infty } F(z) = \lim_{z\rightarrow \infty}
{z^{l-m}\over z^{-m}e^{x} - e^{-x}} =0. \label{zero2}
\end{equation}
This means that for $0<l<m$, $F(z)$ has zeroes at $0$ and at
$\infty$. If the inequalities are non-strict $0\leq l \leq m$,
then the limits are bounded. Similarly for $m\leq 0$, we get that
for $0\geq l \geq m$, $F(z)$ has zeroes at $0$ and at $\infty$,
and if the inequalities are not strict then the limits of $F(z)$
are bounded.

}
\end{remark}

By Remark~\ref{the-remark}, we will have control over the
behaviour of each factor of each summand in $\chii(X,\Ol(\ext^p
D^*\otimes L^{k}))_z$ at $0$ and $\infty$  if
\[
\left|m_{i_1}+\ldots+m_{i_p}+kh\right|
\leq \sum_{i=1}^{2n+1}\vert m_i\vert
\]
for every $p$-tuple $1\leq i_1<\ldots<i_p\leq 2n$.

We will consider three cases which show the general pattern, where $n$
will be as large as needed to illustrate the procedure.

\begin{itemize}
\item \underline{Case $p=0$}.
In this case, we need to determine bounds for $k$ such that
\[
\left|kh\right|
\leq \sum_{i=1}^{2n+1}\vert m_i\vert.
\]
\begin{itemize}
\item When $k=0$ there is nothing to do.
\item When $k=1$ and there exists $i$ such that $\sigma(i)\not = i$, by
\rf{main-identity}
\[
|h|\leq |m_i+m_{\sigma(i)}|\leq |m_i|+|m_{\sigma(i)}|
\leq \sum_{i=1}^{2n+1}|m_i|.
\]
If there is no such $i$, $m_1=m_2=h/2$ so that
\[
|h|= |m_1+m_2|\leq |m_1|+|m_2| \leq \sum_{i=1}^{2n+1}|m_i|.
\]
\item When $k=2$ and there exist $i\not = j$ such that
$\sigma(i)\not = i$ and $\sigma(j)\not = j$, by
\rf{main-identity}
\begin{eqnarray}
|2h|&\leq& |m_i+m_{\sigma(i)}+m_j+m_{\sigma(j)}|\nonumber\\
&\leq& |m_i|+|m_{\sigma(i)}|+|m_j|+|m_{\sigma(j)}|\nonumber\\
&\leq& \sum_{i=1}^{2n+1}|m_i|.\nonumber
\end{eqnarray}
If there is only one $i$ such that $\sigma(i)\not = i$,
there must be $j,k\not =i$ with $m_j=m_k=h/2$ so that
\[
|2h|\leq |m_i+m_{\sigma(i)}+m_j+m_k|
\leq |m_i|+|m_{\sigma(i)}|+|m_j|+|m_k|
\leq \sum_{i=1}^{2n+1}|m_i|.
\]
If there is no such $i$, $m_1=m_2=m_3=m_4=h/2$ so that
\[
|2h|\leq |m_1+m_2+m_3+m_4|
\leq |m_1|+|m_2|+|m_3|+|m_4|
\leq \sum_{i=1}^{2n+1}|m_i|.
\]
\end{itemize}
We continue like this until
$k=n+1$, where the last $h$ is replaced with
$h=m_{2n+1}$ so that
\[0\leq k\leq n+1.\]

\item \underline{Case $p=1$}.
For the sake of simplicity, let us determine
bounds for $k$ such that
\[
\left|m_1+kh\right|
\leq \sum_{i=1}^{2n+1}\vert m_i\vert.
\]
The argument will be analogous for all other $m_i$.
\begin{itemize}
\item When $k=0$ there is nothing to do.
\item When $k=1$ and there exists $i\not =1$ such that $\sigma(i)\not = i$, by
\rf{main-identity}
\[
|m_1+h|\leq |m_1+m_i+m_{\sigma(i)}|\leq |m_1|+|m_i|+|m_{\sigma(i)}|
\leq \sum_{i=1}^{2n+1}|m_i|.
\]
If there is no such $i$, $m_2=m_3=h/2$ so that
\[
|m_1+h|\leq |m_1+m_2+m_3|\leq |m_1|+|m_2|+|m_3|
\leq \sum_{i=1}^{2n+1}|m_i|.
\]
\item When $k=2$ and there exist $i\not = j$ different from 1 and such that
$\sigma(i)\not = i$ and $\sigma(j)\not = j$, by
\rf{main-identity}
\begin{eqnarray}
|m_1+2h|&\leq& |m_1+m_i+m_{\sigma(i)}+m_j+m_{\sigma(j)}|\nonumber\\
&\leq& |m_1|+|m_i|+|m_{\sigma(i)}|+|m_j|+|m_{\sigma(j)}|\nonumber\\
&\leq& \sum_{i=1}^{2n+1}|m_i|.\nonumber
\end{eqnarray}
If there is only one $i$ such that $\sigma(i)\not = i$,
there must be $j,k\not =i,1$ with $m_j=m_k=h/2$ so that
\begin{eqnarray}
|m_1+2h|&\leq& |m_1+m_i+m_{\sigma(i)}+m_j+m_k|\nonumber\\
&\leq& |m_1|+|m_i|+|m_{\sigma(i)}|+|m_j|+|m_k|\nonumber\\
&\leq& \sum_{i=1}^{2n+1}|m_i|.\nonumber
\end{eqnarray}
If there is no such $i\not= 1$, $m_2=m_3=m_4=m_5=h/2$ so that
\begin{eqnarray}
|m_1+2h|&\leq& |m_1+m_2+m_3+m_4+m_5|\nonumber\\
&\leq& |m_1|+|m_2|+|m_3|+|m_4|+|m_5|\nonumber\\
&\leq& \sum_{i=1}^{2n+1}|m_i|.\nonumber
\end{eqnarray}
\end{itemize}
We continue like this until
$k=n$, where the last $h$ is replaced with
$h=m_{2n+1}$ so that
\[0\leq k\leq n.\]

\item \underline{Case $p=n+1$}.
Just as before, let us determine bounds for $k$ such that
\[
\left|m_1+\ldots+m_n+m_{n+1}+kh\right|
\leq \sum_{i=1}^{2n+1}\vert m_i\vert.
\]
The argument will be analogous for all other $(n+1)$-tuples.
\begin{itemize}
\item When $k=0$ there is nothing to do.
\item When $k=1$, in the worst case scenario $\sigma(i)\not = 1,2,\ldots, n$ for
$i= 1,2,\ldots, n$, so that by \rf{h=m2n+1} we can only substitute one $h=m_{2n+1}$.
\begin{eqnarray}
|m_1+\ldots+m_n+m_{n+1}+h|&\leq& |m_1+\ldots+m_n+m_{n+1}+m_{2n+1}| \nonumber\\
&\leq& |m_1|+\ldots+|m_n|+|m_{n+1}|+|m_{2n+1}|\nonumber\\
&\leq& \sum_{i=1}^{2n+1}|m_i|.\nonumber
\end{eqnarray}
\item This time it could also happen that, for instance,
$\sigma(1)=n+1$, and one could subtract one $h=m_1+m_{\sigma(1)}$.
\begin{eqnarray}
|m_1+\ldots+m_n+m_{n+1}-h|&\leq& |m_2+\ldots+m_n| \nonumber\\
&\leq& |m_2|+\ldots+|m_n|\nonumber\\
&\leq& \sum_{i=1}^{2n+1}|m_i|.\nonumber
\end{eqnarray}
\end{itemize}
Thus
\[-1\leq k\leq 1.\]
\end{itemize}
In this fashion, we can obtain all the inequalities stated in the theorem.

\vspace{.1in}

On the one hand, the right hand side in \rf{localization-formula} can be
considered as a meromorphic function with possibly a finite number of poles on the
unit circle  and at the origin.
On the other hand, since $\chii(X,\Ol(\ext^p D^*\otimes L^{-k}  ))_z$ is an index, it is
also a finite Laurent polynomial in $z$ and can be regarded as a meromorphic funtion of the form $\sum_{j} a_{j}z^{j}$, $a_j\in\mathbb{Z}$, for finitely many $j\in\mathbb{Z}$.
Taking the limits
at $0$ and $\infty$ of both sides tells us that $a_{j} =0$ for $j\not=0$ if
the inequalities of Theorem~\ref{main-theorem} are fulfilled. Furthermore, if one side of
the corresponding inequality is strict, $a_0=0$ as well.
\qd

\section{Special circle actions}\label{special-case}

If the complex contact manifold admits a circle action whose
exponents $\{m_i\}$ are all non-negative, then one can prove further
rigidity and vanishing results, such as the following.

\begin{prop}\label{proposition-special-case}
Let $X$ be a complex contact manifold, $D$ the contact
distribution and $L=TX/D$. Assume $X$ admits a circle
action by holomorphic automorphisms preserving the contact
structure, whose exponents $\{m_i\}$
are all non-negative at any $S^1$-fixed point component.
Then, the equivariat holomorphic Euler characteristic
$\chii(X,\Ol(\sym^pD^*\otimes L^{-k}))_z$ is rigid, i.e.
\[\chii(X,\Ol(\sym^pD^*\otimes L^{-k}))_z=\chii(X,\Ol(\sym^pD^*\otimes L^{-k})),\]
$z\in S^1$, if
$0\leq k\leq n+1-p$, for $0\leq p\leq n$,
where $\sym^pD^*$ denotes the $p$-th symmetric tensor power of the
bundle $D^*$.
Furthermore
\[\chii(X,\Ol(\sym^pD^*\otimes L^{-k}))=0\]
if one side of the corresponding inequality is strict.
\end{prop}

{\em Proof}. For the sake of simplicity we will consider the case
$p=2$.

Recall that
\[
m_{2n+1}={1\over n+1}\left( \sum_{i=1}^{2n+1}
m_i\,\,\,\right)=h.
\]
Since all the exponents are non-negative, the
relevant inequalities now take the form
\[
0\leq m_i+m_{j} + k m_{2n+1} \leq (n+1)m_{2n+1},
\]
for $1\leq i\leq  j\leq 2n$.
\begin{itemize}
\item If $\sigma(i) = i$ and $\sigma(j) = j$, the inequality becomes
\[
0\leq (k+1) m_{2n+1} \leq (n+1)m_{2n+1}.
\]
so that $-1\leq k\leq n$.
\item If $\sigma(i) = i$ and $\sigma(j) \not= j$, by \rf{main-identity} $m_i=m_{2n+1}/2$,
$m_j<m_{2n+1}$  and
\[
0\leq m_i + m_{j} + k m_{2n+1} < m_{2n+1} + (k+1/2) m_{2n+1}\leq(n+1)m_{2n+1}.
\]
which requires $0\leq k<n-1/2$.
\item If $\sigma(i) \not= i$ and $\sigma(j) \not= j$, by \rf{main-identity} $m_i<m_{2n+1}$,
$m_j<m_{2n+1}$  and
\[
0\leq m_i + m_{j} + k m_{2n+1} < 2m_{2n+1} + k m_{2n+1}\leq(n+1)m_{2n+1},
\]
which requires $0\leq k\leq n-1$.
\end{itemize}
Therefore $k$ must satisfy
\[0\leq k\leq n-1.\]

Similarly for other values of $p$.
\qd

\begin{remark}{\rm
The homogeneous complex contact manifold
\[Z={SO(8)\over SO(4)\times U(2)}\]
does not admit a holomorphic action with non-negative exponents such as the
one in Proposition~\ref{proposition-special-case}, since by he
Bott-Borel-Weil theorem
\[\chii(Z,\Ol(\sym^2D^*\otimes L^{-1}))=1.\]
}
\end{remark}


{\small
\renewcommand{\baselinestretch}{0.5}
\newcommand{\bi}{\vspace{-.05in}\bibitem} }

\end{document}